\pdfoutput=1
\RequirePackage{ifpdf}
\ifpdf % We are running pdfTeX in pdf mode
\documentclass[pdftex]{sigma}
\else
\documentclass{sigma}
\fi

\usepackage{tikz-cd}
\usepackage{bm}

\newtheorem{thm}{Theorem}[section]

\newtheorem{conj}[thm]{Conjecture}

\newtheorem{prop}[thm]{Proposition}

\theoremstyle{definition}
\newtheorem{ex}[thm]{Example}
\newtheorem{fact}[thm]{Fact}
\newtheorem{rem}[thm]{Remark}
\newtheorem{defn}[thm]{Definition}
\newtheorem{exercise}[thm]{Exercise}
\newtheorem{defnthm}[thm]{Definition/Theorem}

\numberwithin{equation}{section}

\newcommand{\cB}{\ensuremath{\mathcal B}}
\newcommand{\cL}{\ensuremath{\mathcal L}}

\newcommand{\cF}{\ensuremath{\mathcal F}}

\newcommand{\cM}{\ensuremath{\mathcal M}}
\newcommand{\cO}{\ensuremath{\mathcal O}}

\newcommand{\R}{\ensuremath{\mathbb R}}
\newcommand{\C}{\ensuremath{\mathbb C}}
\newcommand{\CP}{\ensuremath{\mathbb {CP}}}

\newcommand{\Z}{\ensuremath{\mathbb Z}}

\newcommand{\GL}{\mathrm{GL}}

\newcommand{\SU}{\mathrm{SU}}
\newcommand{\U}{\mathrm{U}}

\newcommand{\cE}{{\mathcal E}}

\newcommand{\delbar}{\ensuremath{\overline{\partial}}}
\newcommand{\zbar}{\ensuremath{\overline{z}}}

\newcommand{\semif}{\ensuremath{\mathrm{sf}}}
\newcommand{\app}{\ensuremath{\mathrm{app}}}

\newcommand{\I}{{\mathrm i}}
\newcommand{\e}{{\mathrm e}}
\newcommand{\de}{\mathrm{d}}

\newcommand{\norm}[1]{\lVert#1\rVert}
\newcommand{\IP}[1]{\langle#1\rangle}

\newcommand{\eps}{\epsilon}

\newcommand{\del}{{\partial}}

\DeclareMathOperator{\Det}{Det}
\DeclareMathOperator{\End}{End}

\DeclareMathOperator{\Aut}{Aut}

\begin{document}

\allowdisplaybreaks

\newcommand{\arXivNumber}{1809.05735}

\renewcommand{\thefootnote}{}

\renewcommand{\PaperNumber}{018}

\FirstPageHeading

\ShortArticleName{Perspectives on the Asymptotic Geometry of the Hitchin Moduli Space}

\ArticleName{Perspectives on the Asymptotic Geometry\\ of the Hitchin Moduli Space\footnote{This paper is a~contribution to the Special Issue on Geometry and Physics of Hitchin Systems. The full collection is available at \href{https://www.emis.de/journals/SIGMA/hitchin-systems.html}{https://www.emis.de/journals/SIGMA/hitchin-systems.html}}}

\Author{Laura FREDRICKSON}

\AuthorNameForHeading{L.~Fredrickson}

\Address{Stanford University, Department of Mathematics, 380 Serra Mall, Stanford, CA 94305, USA}
\Email{\href{mailto:lfredrickson@stanford.edu}{lfredrickson@stanford.edu}}
\URLaddress{\url{https://web.stanford.edu/~ljfred4/}}

\ArticleDates{Received September 23, 2018, in final form February 25, 2019; Published online March 11, 2019}

\Abstract{We survey some recent developments in the asymptotic geometry of the Hitchin moduli space, starting with an introduction to the Hitchin moduli space and hyperk\"ahler geometry.}

\Keywords{Hitchin moduli space; Higgs bundles; hyperk\"ahler metric}

\Classification{53C07; 53C26}

\renewcommand{\thefootnote}{\arabic{footnote}}
\setcounter{footnote}{0}

\medskip

Fix a compact Riemann surface $C$. In his seminal paper ``The self-duality equations on a~Riemann surface''~\cite{Hitchin87}, Hitchin introduced the moduli space $\cM$ of ${\rm SL}(2,\C)$-Higgs bundles on~$C$ and proved that $\cM$ admits a hyperk\"ahler metric~$g_{\cM}$. In these notes, we give an introduction to the hyperk\"ahler geometry of the Hitchin moduli space, focusing on the geometry of the ends of the Hitchin moduli space. In the last section (Section~\ref{sec:4}), we briefly survey some recent developments in the description of the asymptotic geometry of $\cM$. We start with Gaiotto--Moore--Neitzke's conjectural description in~\cite{GMNwallcrossing, GMNhitchin} and survey recent progress in \cite{DumasNeitzke, Fredricksonasymptoticgeometry, FredricksonSLn, MSWW14, MSWW17}. We take a meandering path through more classical geometric results to get there. In Section~\ref{sec:1}, we give a survey of the results in Hitchin's original paper~\cite{Hitchin87}, since many current lines of research originate there. In Section~\ref{sec:2}, we focus on the hyperk\"ahler metric on the Hitchin moduli space. In order to more fully describe the conjectured picture of the hyperk\"ahler metric on the Hitchin moduli space, we take a detour into the classification of noncompact hyperk\"ahler 4-manifolds as ALE, ALF, ALG, and ALH, highlighting some classical and more recent results. In Section~\ref{sec:3}, we consider the spectral data of the Hitchin moduli space. We describe the abelianization of Hitchin's equations near the ends of the moduli space, and the resulting importance of the spectral data for the asymptotic geometry.

\section{A tour of the Hitchin moduli space}\label{sec:1}

Given the data of
\begin{itemize}\itemsep=0pt
 \item $C$, a compact Riemann surface of genus $\gamma_C \geq 2$ (unless indicated otherwise), and
 \item $E \rightarrow C$, a complex vector bundle of rank $n$
\end{itemize}
we get a Hitchin moduli space $\cM$.

In one avatar, the Hitchin moduli space is the moduli space of $\GL(n,\C)$-Higgs bundles up to equivalence. In another avatar, the Hitchin moduli space is the moduli space of $\GL(n,\C)$-flat connections up to equivalence. In this section, we define these objects and explain the correspondence between Higgs bundles (a holomorphic object) and flat connections and their associated representations (a representation theoretic object). We survey many of the results appearing in Nigel Hitchin's seminal paper for the complex Lie group $G_\C={\rm SL}(2,\C)$~\cite{Hitchin87}.

\subsection{Motivation: Narasimhan--Seshadri correspondence}

For the sake of motivation, there is an earlier example of a correspondence between holomorphic objects and representations. In 1965, Narasimhan and Seshadri proved the equivalence between stable holomorphic vector bundles
on a compact Riemann surface $C$ and irreducible projective unitary representations of the fundamental group~\cite{NarasimhanSeshadri}. In 1983, Donaldson gave a more direct proof of this fact using the differential geometry of connections on holomorphic bundles~\cite{DonaldsonNS}. We specialize to the degree $0$ case for simplicity, so that projective unitary representations are simply unitary representations.

\begin{thm}[\cite{DonaldsonNS}]\label{thm:DNS} Let $\cE$ be a indecomposable holomorphic bundle of $\deg \cE=0$ and rank $n$ over a Riemann surface~$C$. The holomorphic bundle $\cE$ is stable if, and only if, there is an irreducible flat unitary connection on $\cE$. Taking the holonomy representation of this connection, we have the following equivalence:
\begin{align*}
 \left\{\begin{matrix} \mbox{stable holomorphic bundles}\\\cE \end{matrix} \right\}_{\!\mbox{$\!/\!\!\sim$}} &\leftrightarrow
 \left\{\begin{matrix} \mbox{flat $\mathrm{U}(n)$-connections}\\\nabla \end{matrix} \right\}_{\!\mbox{$\!/\!\!\sim$}} \\
 & \leftrightarrow
 \left\{\begin{matrix} \mbox{irreducible representations}\\ \rho\colon \pi_1(C) \rightarrow \mathrm{U}(n) \end{matrix} \right\}_{\!\mbox{$\!/\!\!\sim$}}.
\end{align*}
\end{thm}
This map from a holomorphic bundle $\cE$ to a flat connection $\nabla$ features a distinguished hermitian metric $h$ on $\cE$. First, note that given any hermitian metric $h$ on a holomorphic bundle $\cE$, there is a unique connection $D(\delbar_\cE, h)$ called the Chern connection characterized by the property that (1) $D^{0,1}=\delbar_\cE$ and (2) $D$ is unitary with the respect to~$h$, i.e., $\de \IP{s_1, s_2}_h=\IP{Ds_1, s_2}_h + \IP{s_1, Ds_2}_h$. The proof of Theorem~\ref{thm:DNS} relies on the following fact: given a stable holomorphic bundle~$\cE$ of degree~$0$, there is a hermitian metric~$h$~-- known as a Hermitian--Einstein metric~-- such that the Chern connection is flat. Consequently, the flat connection associated to~$\cE$ is $\nabla=D(\delbar_E, h)$, where $h$ is the Hermitian--Einstein metric.

The nonabelian Hodge correspondence, which interpolates between the two avatars of the Hitchin moduli space, also features a distinguished hermitian metric.

\subsection{Definition of the Hitchin moduli space} Higgs bundles and the Hitchin moduli space first appeared in Nigel Hitchin's beautiful paper ``The self-duality equations on a Riemann surface''~\cite{Hitchin87}. We specialize to the degree $0$ case for simplicity.

\begin{defn}Fix a complex vector bundle $E \rightarrow C$ of degree $0$. A \emph{Higgs bundle} on $E \rightarrow C$ is a pair $(\delbar_E, \varphi)$ where
\begin{itemize}\itemsep=0pt
 \item $\delbar_E$ is a holomorphic structure on $E$ (We'll denote the corresponding holomorphic vector bundle by $\cE =(E, \delbar_E)$.)
 \item $\varphi \in \Omega^{1,0}(C, \operatorname{End} E)$ is called the ``Higgs field''
\end{itemize}
satisfying $\delbar_E \varphi =0$. (Alternatively, the Higgs field is a holomorphic map $\varphi\colon \cE \rightarrow \cE \otimes K_C$, where $K_C=\mathcal{T}^{1,0}(C)$ is the canonical bundle.)
\end{defn}
\begin{rem} If we want $G_\C={\rm SL}(n,\C)$ rather than $\GL(n,\C)$, we must impose the condition $\operatorname{tr} \varphi =0$, since $\mathfrak{sl}(n,\C)$ consists of traceless matrices. Additionally, we insist that $\Det \cE \simeq \cO_C$, as holomorphic line bundles.
\end{rem}

\begin{defn} Fix a Higgs bundle $(\cE, \varphi)$. A hermitian metric on $E$, the underlying complex vector bundle, is \emph{harmonic} if
 \begin{gather*}
 F_{D(\delbar_E, h)} + \big[\varphi, \varphi^{*_h}\big]=0.
 \end{gather*}
 Here $F_D$ is the curvature of $D$; $\varphi^{*_h}$ is the hermitian adjoint\footnote{In a local holomorphic coordinate $z$ and a local holomorphic frame for $(E, \delbar_E)$, if $\varphi = \Phi \de z$, then $\varphi^{*_h} = h^{-1} \Phi^* h \de \zbar$.} of $\varphi$ with respect to~$h$.
\end{defn}

\begin{defn}A triple $(\delbar_E, \varphi, h)$ is a \emph{solution of Hitchin's equations} if $(\delbar_E, \varphi)$ is a Higgs bundle and $h$ is harmonic, i.e.,
 \begin{gather} \label{eq:Hitchineq}
 \delbar_E \varphi =0, \qquad F_{D(\delbar_E, h)} + \big[\varphi, \varphi^{*_h}\big]=0.
 \end{gather}
\end{defn}

\begin{defn} Fix a complex vector bundle $E \rightarrow C$. The associated \emph{Hitchin moduli space} $\cM$ consists of triples $\big(\delbar_E, \varphi, h\big)$ solving Hitchin's equations, up to complex gauge equivalence
 \begin{gather*}
 g \cdot (\delbar_E, \varphi, h) = \big(g^{-1} \circ \delbar_E \circ g, g^{-1} \varphi g, g \cdot h \big), \qquad \mbox{where $(g \cdot h)(v,w) = h(g v, gw)$},
 \end{gather*}
 for $g \in \Gamma(C, \Aut(E))$.
\end{defn}

The Hitchin moduli space is a manifold with singularities. When $\gamma_C \geq 2$, the dimension of the $\U(n)$-Hitchin moduli space is $\dim_\R \cM(C, \U(n))=4\big(n^2(g-1) + 1\big)$; the dimension of the $\SU(n)$-Hitchin
moduli space is $\dim_\R \cM(C, \SU(n))=4\big(n^2-1\big)(g-1)$.

\begin{exercise} \label{ex:model}Verify that the following triple $(\delbar_E, t\varphi, h_t)$ on $\C$ solves Hitchin's equations:
\begin{gather*} \delbar_E =\delbar, \qquad t\varphi = t\begin{pmatrix} 0 & 1 \\ z & 0 \end{pmatrix} \de z, \qquad h_t=\begin{pmatrix} |z|^{1/2} \e^{u_t(|z|)} &\\
 & |z|^{-1/2} \e^{-u_t(|z|)}
 \end{pmatrix}, \end{gather*}
 where $u_t=u_t(|z|)$ is the solution of the ODE
 \begin{gather*}
 \left(\frac{\de^2}{\de |z|^2} + \frac{1}{|z|}\frac{\de}{\de |z|} \right) u_t= 8 t^2 |z|\sinh(2 u_t),
 \end{gather*}
with boundary conditions $u_t(|z|) \sim -\frac{1}{2} \log |z|$ near $|z|=0$ and $\lim\limits_{|z| \to \infty} u_t(|z|)=0$. It may be useful to note:
\begin{itemize}\itemsep=0pt
 \item In a local holomorphic frame where $\delbar_E=\delbar$, the curvature is $F_{D(\delbar_E, h)}= \delbar \big(h^{-1} \del h \big)$. When~$h$ is diagonal, $F_{D(\delbar_E, h)}=\delbar \del \log h$.
 \item Let $z=x+{\rm i}y$ be a local holomorphic coordinate. Then, $\delbar \del \nu = \frac{1}{4} \big(\frac{\de^2}{\de x^2} + \frac{\de^2}{\de y^2} \big) \nu \, \de \zbar \wedge \de z$.
\end{itemize}

{\bf Note:} This is the model solution featured in \cite{FredricksonSLn, GMNhitchin, MSWW14}. The base curve is $\CP^1$ with an irregular singularity at~$\infty$~\cite{FredricksonNeitzke}.
\end{exercise}

\subsection{Nonabelian Hodge correspondence}

The Hitchin moduli space $\cM$ is hyperk\"ahler. As a consequence, it has a $\CP^1$-worth of complex structures, labeled by parameter $\zeta \in \CP^1$. Two avatars of the Hitchin moduli space are
\begin{itemize}\itemsep=0pt
 \item the Higgs bundle moduli space ($\zeta=0$), and
 \item the moduli space of flat $\GL(n,\C)$-connections $\zeta \in \C^\times$.
\end{itemize}

Starting with the triple $[(\delbar_E, \varphi, h)]$ in $\cM$, the associated Higgs bundle $\big[\big(\delbar_E, \varphi\big)\big]$ is obtained by forgetting the harmonic metric~$h$. Starting with the triple $\big[\big(\delbar_E, \varphi, h\big)\big]$, for each $\zeta \in \C^\times$, the associated flat connection is $[\nabla_\zeta]$ where
\begin{gather}\label{eq:nabla}
 \nabla_\zeta = \zeta \varphi + D_{(\delbar_E, h)} + \zeta^{-1} \varphi^{*_h}.
\end{gather}
The nonabelian Hodge correspondence describes the correspondence between solutions of Hit\-chin's equations, Higgs bundles, and flat connections. It answers questions that include ``What Higgs bundles admit harmonic metrics?'' and ``Can any flat connection be produced in this way?''

\begin{exercise}Use Hitchin's equations in \eqref{eq:Hitchineq} to verify that $\nabla_\zeta$ in \eqref{eq:nabla} is flat.
\end{exercise}

What Higgs bundles $(\cE, \varphi)$ admit harmonic metrics $h$? The following algebraic stability condition guarantees the \emph{existence} of a harmonic metric. Moreover, any harmonic metric on an indecomposable Higgs bundle is \emph{unique} up to rescaling by a constant. We define stability for holomorphic bundles, before generalizing it to Higgs bundles.

\begin{defn} A holomorphic bundle $\cE$ is \emph{stable} if for every proper holomorphic subbundle $\mathcal{F} \subset \cE$, the slopes $\mu(\mathcal{F}):= \frac{\deg \cF}{\operatorname{rank} \cF}$
 satisfy
 \begin{gather*}
 \mu(\mathcal{F}) < \mu(\mathcal{E}).
 \end{gather*}
\end{defn}

\begin{defn}A Higgs bundle $(\cE, \varphi)$ is \emph{stable} if for every $\varphi$-invariant proper holomorphic subbundle $\mathcal{F} \subset \cE$, the slopes satisfy
 \begin{gather*}
 \mu(\mathcal{F}) < \mu(\mathcal{E}).
 \end{gather*}
A Higgs bundle $(\cE, \varphi)$ is \emph{polystable} if it is the direct sum of stable Higgs bundles of the same slope.
\end{defn}

\begin{thm}[\cite{Hitchin87, simpsonhiggs}] A Higgs bundle admits a harmonic metric if, and only if, it is poly\-stable.
\end{thm}

The nonabelian Hodge correspondence gives an equivalence between Higgs bundles, solutions of Hitchin's equations, and flat connections. Admittedly, our presentation in this paper focuses on the equivalence between Higgs bundles and solutions of Hitchin's equations, while neglecting flat connections. For more on the equivalence between flat connections and solutions of Hitchin's equations, see, for example,~\cite{Wentworth2014}.

\begin{thm}[nonabelian Hodge correspondence, \cite{Corlette, donaldsonSD, Hitchin87, simpsonhiggs}] Fix a complex vector bundle $E \rightarrow C$ of rank $n$ and degree $0$, and take $G_\C={\rm SL}(n,\C)$. There is a correspondence between polystable ${\rm SL}(n,\C)$-Higgs bundles and completely reducible ${\rm SL}(n,\C)$-connections\footnote{Let $E$ be a complex vector bundle. A connection $\nabla$ is called \emph{completely reducible} if every $\nabla$-invariant subbundle $F \subset E$ has a $\nabla$-invariant complement. A connection $\nabla$ is called \emph{irreducible} if there are no nontrivial proper $\nabla$-invariant subbundles.}:
\begin{align*}
 \left\{\begin{matrix} \mbox{polystable Higgs bundle}\\ \big(\delbar_E, \varphi\big) \end{matrix} \right\}_{\!\mbox{$\!/\!\!\sim$}}
& \leftrightarrow \left\{\begin{matrix} \mbox{soln of Hitchin's eq}\\ \big(\delbar_E, \varphi, h\big) \end{matrix} \right\}_{\!\mbox{$\!/\!\!\sim$}}\\
& \leftrightarrow
 \left\{\begin{matrix} \mbox{completely reducible} \\ \mbox{flat ${\rm SL}(n,\C)$-connection $\nabla$} \end{matrix} \right\}_{\!\mbox{$\!/\!\!\sim$}}.
\end{align*}
In this correspondence $\big(\delbar_E, \varphi\big)$ is stable if, and only if, the associated flat connection is irreducible; this is the smooth locus of $\cM$.
\end{thm}

In the above correspondence, we typically associate the connection $\nabla_{\zeta=1}$ from~\eqref{eq:nabla}. To get a representation, we use the Riemann--Hilbert equivalence
\begin{gather*}
 \left\{\begin{matrix} \mbox{flat ${\rm SL}(n,\C)$-connections}\\\nabla \end{matrix} \right\}_{\!\mbox{$\!/\!\!\sim$}} \longleftrightarrow
 \left\{\begin{matrix} \mbox{representation}\\ \rho\colon \pi_1(C) \rightarrow {\rm SL}(n,\C) \end{matrix} \right\}_{\!\mbox{$\!/\!\!\sim$}}.
\end{gather*}
To go from a flat connection to a representation, simply take the monodromy of a connection. In the other direction, to go from a representation~$\rho$ to a bundle with flat connection, take the trivial bundle $\underline{\C^n} \rightarrow \widetilde{C}$ on the universal cover $\widehat{\pi}\colon \widetilde{C} \to C$ and equip it with the trivial flat connection given by exterior differentiation. The bundle with flat connection on $C$ is obtained by quotienting by the following equivalence relation on pairs $(x, v) \in \widetilde{C} \times \C^n$: for any $\gamma \in \pi_1(C)$,
\begin{gather*}
 (x, v) \sim \big(\widehat{\pi}^*\gamma \cdot x, \rho(\gamma) v\big).
\end{gather*}
Here, $\widehat{\pi}^*\gamma$ is the path in $\widetilde{C}$; it is the lift of $\gamma$ with initial point $x\in \widetilde{C}$; $\widehat{\pi}^*\gamma \cdot x$ is the terminal point of the path $\pi^*\gamma$.

\begin{exercise} Describe the Higgs bundles $(\cE, \varphi)$ in the $\GL(1,\C)$-Higgs bundle moduli space over $C$.
\end{exercise}

\begin{exercise}\quad
\begin{itemize}\itemsep=0pt
 \item[(a)] Describe $\chi_{{\rm SL}(2,\C)}\big(T^2\big)$, the ${\rm SL}(2,\C)$ character variety of~$T^2$.
 \item[(b)] Describe an isomorphism $\psi\colon \big(\C^\times \times \C^\times\big) /\sigma\rightarrow \chi_{{\rm SL}(2,\C)}\big(T^2\big)$ where $\sigma\colon (a,b)\mapsto(-a,-b)$.
\end{itemize}
\end{exercise}

\subsection{Hitchin fibration} \label{sec:spectraldata}
\begin{figure}[ht]\centering
\includegraphics[height=1.0in]{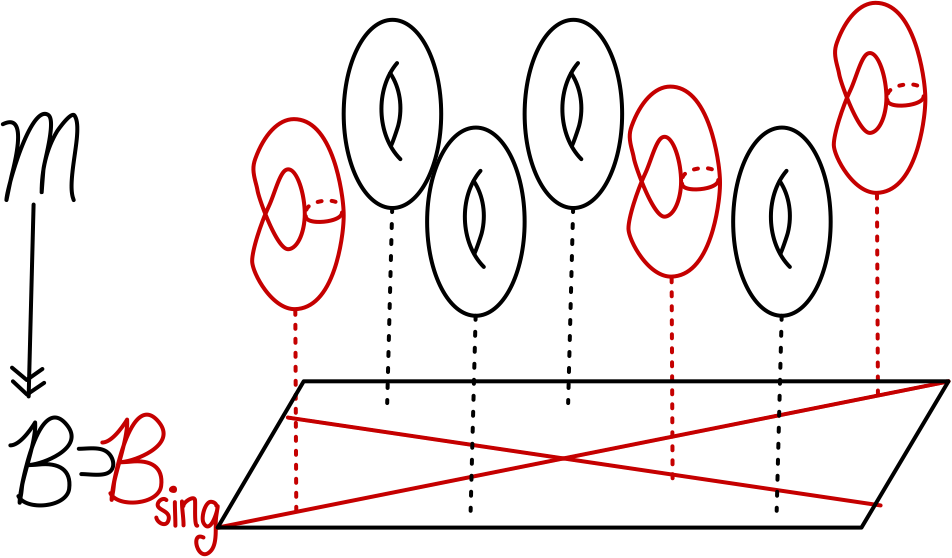}
\caption{ Hitchin fibration.}\label{fig:Hitchinfibration}
\end{figure}

The Hitchin fibration is a surjective holomorphic map
\begin{align} \label{eq:Hitchinfibration}
 \operatorname{Hit}\colon \ \cM &\twoheadrightarrow \cB \simeq \C^{\frac{1}{2} \dim_\C \cM},\\ \nonumber
 (\delbar_E, \varphi, h) &\mapsto \operatorname{char}_\varphi (\lambda) \quad \mbox{[encodes eigenvalues of $\varphi$]},
\end{align}
where $\operatorname{char}_\varphi(\lambda)$ is the characteristic polynomial of $\varphi$. Fundamentally, the Hitchin fibration $\operatorname{Hit}$ maps the Higgs field $\varphi$ to its eigenvalues $\lambda_1, \dots, \lambda_n$ (multivalued sections of~$K_C)$. With the map~$\operatorname{Hit}$, $\cM$~is a~``an algebraic completely integrable system''\footnote{An algebraic completely integrable system~$\cM$ is a holomorphic symplectic space fibered over a complex base~$\cB$ with $\dim_\C \cB = \frac{1}{2} \dim_C \cM$; the fibers are Lagrangian; generic fibers are abelian varieties~\cite{Donagi:1995am}.}. The $\frac{1}{2}$-dimensional compact complex torus fibers degenerate over a complex codimension-one locus $\cB_{{\rm sing}}$, as indicated in Fig.~\ref{fig:Hitchinfibration}. The most singular fiber, $\operatorname{Hit}^{-1}(\mathbf{0}) \subset \cM$, is called the ``nilpotent cone'', and it contains the space of stable holomorphic vector bundles. Let $\cB'=\cB-\cB_{{\rm sing}}$ and call $\cM'=\operatorname{Hit}^{-1}(\cB')$ the \emph{regular locus} of the Hitchin moduli space. It is obvious that the Hitchin moduli space $\cM$ is noncompact, since $\cB$ is noncompact.

Specializing to the case ${\rm SL}(2,\C)$, note that
\begin{gather*}
\operatorname{char}_\varphi (\lambda) = (\lambda-\lambda_1)(\lambda-\lambda_2) = \lambda^2 -(\lambda_1+\lambda_2) \lambda + \lambda_1\lambda_2 = \lambda^2 - \operatorname{tr} \varphi \lambda + \det \varphi.
\end{gather*}
Note that $\operatorname{tr} \varphi =0$ and $\det \varphi \in H^0\big(C, K_C^2\big)$. Consequently, the Hitchin base $\cB$ is parametrized by the space of holomorphic quadratic differentials $H^0\big(C, K_C^2\big)$.

\begin{exercise}Use the Riemann--Roch formula to verify directly that the complex dimension of
$\cB_{{\rm SL}(2,\C)}$ is $3(g-1)$.

\texttt{Hint:} \emph{The Riemann--Roch formula for line bundles $\cL \rightarrow C$ states that}
\begin{gather*}
 h^0(C,\cL) -h^0\big(C, \cL^{-1} \otimes K_C\big)=\deg(\cL)+1-g,
\end{gather*}
\emph{where $h^0(C,\cL)$ is the dimension of $H^0(C, \cL)$, the space of holomorphic sections of $\cL$. Additio\-nal\-ly, $\deg(K_C)=2g-2$.}
\end{exercise}

The Hitchin fibration has a collection of distinguished sections, known as ``Hitchin sections''. For the ${\rm SL}(2,\C)$-Hitchin moduli space, there are $2^{2\gamma_C}$-Hitchin sections labeled by a choice of a~spin structure on $C$. Given a spin structure $K_C^{1/2}$, the corresponding Hitchin section is
 \begin{gather*}
 \cB \rightarrow \cM, \\
 q_2 \mapsto \cE = K_C^{-1/2} \oplus K_C^{1/2}, \qquad \varphi = \begin{pmatrix} 0 & 1 \\ q_2 & 0 \end{pmatrix}, \qquad h =\begin{pmatrix} h_{K_C^{-1/2}} & 0 \\ 0 & h_{K_C^{1/2}} \end{pmatrix}.
 \end{gather*}
To interpret $\varphi$, view ``$1$'' as the identity map $K_C^{1/2} \rightarrow K_C^{-1/2} \otimes K_C \simeq K_C^{1/2}$, and view tensoring by $q_2$ as a map $K_C^{-1/2} \rightarrow K_C^{1/2} \otimes K_C \simeq K_C^{-1/2} \otimes K_C^2$. The hermitian metric respects this direct sum, and the metric component $h_{K_C^{-1/2}}=h_{K_C^{1/2}}^{-1}$ is determined from Hitchin's equations.

The Hitchin section is related to uniformization. From $h_{K_C^{-1/2}}$, we get a hermitian met\-ric~$h_{K_C^{-1}}$ on the inverse of the holomorphic tangent bundle $K_C^{-1}=\big(\mathcal{T}^{1,0}(C)\big)^{-1}$. Note that this bundle is related to the usual tangent bundle~$TC$. From \cite[Theorem~11.2]{Hitchin87},
 \begin{gather} \label{eq:uniformization}
 g = q_2 + \left(h_{K_C^{-1}} + \frac{|q_2|^2}{h_{K_C^{-1}}}\right) + \overline{q_2}
 \end{gather}
is a Riemannian metric on $C$ of Gaussian curvature $-4$. The map between Teichm\"uller space $\mathrm{Teich}(C)$ and $H^0\big(C, K_C^2\big)$ is further discussed in \cite[Section~3]{WolfTeich}. Note that if $q_2=0$, then the Riemannian metric~$g$ in~\eqref{eq:uniformization} belongs to the conformal class given by the complex structure on~$C$. This is the ``uniformizing metric'' and the corresponding Higgs bundle in~\eqref{eq:uniformizingHiggs} is called the ``uniformizing point''.

\begin{exercise}Consider the Higgs bundle
\begin{gather} \label{eq:uniformizingHiggs}
 \cE = K_C^{-1/2} \oplus K_C^{1/2}, \qquad \varphi = \begin{pmatrix} 0 & 1 \\ 0 & 0 \end{pmatrix},
\end{gather}
where $1$ is the identity map $K_C^{1/2} \to K_C^{-1/2} \otimes K_C$.
\begin{itemize}\itemsep=0pt
\item [(a)] Show that the holomorphic bundle $\cE$ is unstable by exhibiting a destabilizing subbundle, i.e., a holomorphic subbundle $\cL$ such that
\begin{gather*}
 \mu(\cL) \geq \mu(\cE).
\end{gather*}
It might be helpful to note that $\deg K_C=2\gamma_C-2$, where $\gamma_C \geq 2$ is the genus.
\item[(b)] Describe the group of automorphisms of $K_C^{-1/2} \oplus K_C^{1/2}$. Show that the destabilizing bundle from (a) is unique, i.e., it is preserved by all holomorphic automorphisms.
\item[(c)] Show that the Higgs bundle $(\cE, \varphi)$ is stable. Where is the condition ``$\gamma_C \geq 2$'' used?
\end{itemize}
\end{exercise}

\subsection[$\U(1)$-action and topology]{$\boldsymbol{\U(1)}$-action and topology}

There is a $\C^\times$-action on the Higgs bundle moduli space given by
\begin{gather*}
 \xi \in \C^\times\colon \ \big[(\delbar_E, \varphi)\big] \mapsto \big[(\delbar_E, \xi \varphi)\big].
\end{gather*}
(Here $[(\delbar_E, \varphi)]$ denotes the equivalence class in $\cM$.) Similarly, we get a $\U(1)$-action on the Hitchin moduli space:
\begin{gather*}
 \e^{{\rm i} \theta} \in \U(1)\colon \ \big[(\delbar_E, \varphi, h)\big] \mapsto \big[(\delbar_E, \e^{{\rm i} \theta} \varphi, h)\big].
\end{gather*}
The $\U(1)$-action preserves the K\"ahler form on $\cM$ and generates a moment map\footnote{Suppose $(X, \omega)$ is a symplectic manifold with $G$-action. Then for any $Z \in \mathfrak{g}=\operatorname{Lie} G$, we get an associated vector field $\mathfrak{X}_Z$ on~$X$. A function $\mu\colon X \rightarrow \mathfrak{g}^*$ is a moment map for the $G$-action if $\mu$ is $G$-equivariant, and for all $Z \in \mathfrak{g}$, then \begin{gather*}\iota_{\mathfrak{X}_Z} \omega = \de \mu_Z.\end{gather*}

In the case where $G=\U(1)$, $\mathfrak{g}=\I \R$, so (ignoring $\I$) we can view $\mu\colon X \rightarrow \R$ as an ordinary function.}
\begin{gather*}
 \mu = \int_C \operatorname{tr} \big(\varphi \wedge \varphi^{*_h}\big).
\end{gather*}

\begin{figure}[ht]\centering
\includegraphics[height=1.2in]{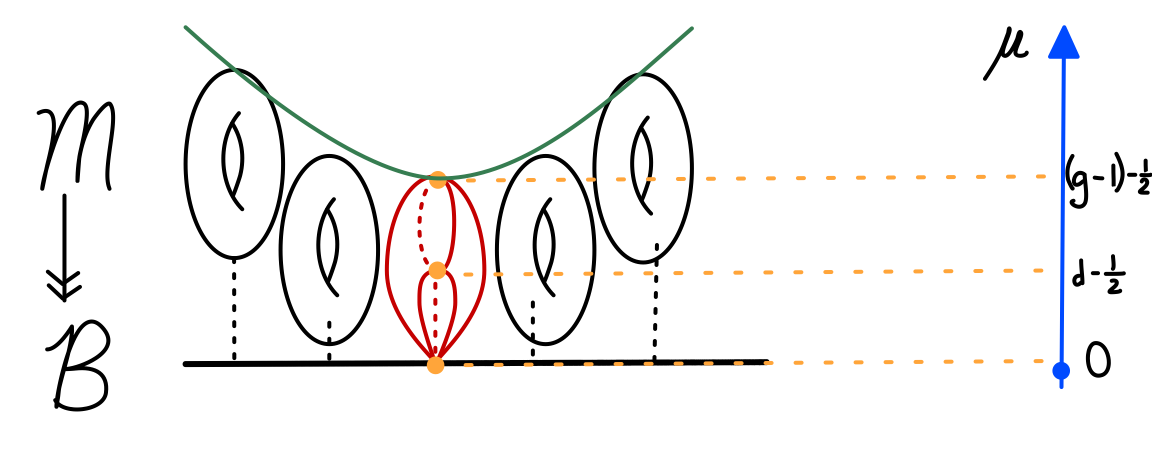}
\caption{The maximal value of $\mu\colon \cM \rightarrow \R$ in each fiber of $\operatorname{Hit}\colon \cM \rightarrow \cB$ is achieved on each of the~$2^{2\gamma_C}$ Hitchin sections.}\label{fig:Morse}
\end{figure}

We specialize to $G_\C = {\rm SL}(2,\C)$ for the rest of Section~\ref{sec:1}. As shown in Fig.~\ref{fig:Morse}, the maximal value of $\mu$ in each torus fiber is achieved on each of the Hitchin sections~\cite{Tholozan}. There are subspaces of $\U(1)$-fixed points. The associated values of $\mu$ are $0$, and $d-\frac{1}{2}$ for $d=1, \dots, \gamma_C-1$, where~$\gamma_C$ is the genus of~$C$. The $\U(1)$-fixed points in $\mu^{-1}(0)$ are the polystable vector bundles. The submanifold of~$\U(1)$-fixed points in $\mu^{-1}\big(d-\frac{1}{2}\big)$ is the space of Higgs bundles
\begin{gather*}
 \cE = \cL^{-1} \oplus \cL, \qquad \varphi = \begin{pmatrix} 0 & \alpha \\ 0 & 0 \end{pmatrix},
\end{gather*}
where $\deg \cL=d$ and $\alpha \in H^0\big(C, \cL^{-2} \otimes K_C\big)$. If $d=\gamma_C-1$, then $\cL = K_C^{1/2}$ and $\mu^{-1}\big(d-\frac{1}{2}\big)$ consists of $2^{2\gamma_C}$ $\U(1)$-fixed points corresponding the $2^{2 \gamma_C}$ choices of spin structure $K_C^{1/2}$ on $C$. Note that each of these $2^{2\gamma_C}$ Higgs bundles described in~\eqref{eq:uniformizingHiggs} gives a different representation of~$\pi_1(C)$ in ${\rm SL}(2,\C)$, however, all project to the same uniformizing representation of $\pi_1(C)$ in $\mathrm{PSL}(2,\C)$.

The topology of the ${\rm SL}(2,\C)$ character variety was originally computed using the $\U(1)$-action on the $\SU(2)$-Hitchin moduli space~$\cM$ in~\cite{Daskalopolousetal}. Since $\cM$ deformation retracts onto the nilpotent cone $\operatorname{Hit}^{-1}(\mathbf{0})$, the topology of $\cM$ is the same as $\operatorname{Hit}^{-1}(\mathbf{0})$. By computing the indices of the $\U(1)$-fixed submanifolds in $\operatorname{Hit}^{-1}(\mathbf{0})$, one can determine the topology of $\operatorname{Hit}^{-1}(\mathbf{0})$ using Morse--Bott theory. (See~\cite{SteveSIGMA} for a more thorough introduction to the topology.)

\subsection[${\rm SL}(2,\R)$-Higgs bundles]{$\boldsymbol{{\rm SL}(2,\R)}$-Higgs bundles}

Recall, the nonabelian Hodge correspondence gives us the following equivalence:
\begin{equation*}
 \left\{\begin{matrix} \mbox{stable ${\rm SL}(2,\C)$-Higgs bundles}\\(\delbar_E, \varphi) \end{matrix} \right\}_{\!\mbox{$\!/\!\!\sim$}} \! \quad \longleftrightarrow \quad
 \left\{\begin{matrix} \mbox{irreducible representations}\\ \rho\colon \pi_1(C) \rightarrow {\rm SL}(2,\C) \end{matrix} \right\}_{\!\mbox{$\!/\!\!\sim$}}.
\end{equation*}
One can define an ${\rm SL}(2,\R)$-Higgs bundle as a ${\rm SL}(2,\C)$-Higgs bundle which correspond to a~${\rm SL}(2,\R)$-representation:
\begin{gather*}
 \overset{\bigcup}{\left\{\begin{matrix} \mbox{${\rm SL}(2,\R)$-Higgs bundles}\\(\delbar_E, \varphi) \end{matrix} \right\}}_{\!\mbox{$\!/\!\!\sim$}} \quad \longleftrightarrow \quad
 \overset{\bigcup}{ \left\{\begin{matrix} \mbox{irreducible representations}\\ \rho\colon \pi_1(C) \rightarrow {\rm SL}(2,\R) \end{matrix} \right\}}_{\!\mbox{$\!/\!\!\sim$}}.
\end{gather*}
The Lie subalgebra $\mathfrak{sl}(2,\R) \subset \mathfrak{sl}(2,\C)$ is preserved by the map $\Phi \rightarrow \overline{\Phi}$. Consequently, ${\rm SL}(2,\R)$-Higgs bundles can be viewed as ${\rm SL}(2,\C)$-Higgs bundles with additional conditions:
\begin{itemize}\itemsep=0pt
 \item $\cE$ has an orthogonal structure $Q\colon \cE \rightarrow \cE^*$, and
 \item $\varphi$ is $Q$-symmetric, i.e., $\varphi^T Q = Q \varphi$,
\begin{equation*}
\begin{tikzcd}
\cE \arrow{r}{\varphi}
\arrow{rd}{Q}
&\cE \otimes K_C \arrow{rd}{Q}\\
& \cE^* \arrow{r}{\varphi^T} & \cE^*\otimes K_C.
\end{tikzcd}
\end{equation*}
\end{itemize}
Note that the harmonic metric $h$ will in turn satisfy $h^T Q h =Q$. We have already encountered some ${\rm SL}(2,\R)$-Higgs bundles. Namely, all Higgs bundles in the Hitchin sections are ${\rm SL}(2,\R)$-Higgs bundles. To see this, just take the orthogonal structure $Q= \left(\begin{smallmatrix} 0 & 1 \\ 1 & 0 \end{smallmatrix}\right)$ where the~``$1$''s represent the identity maps $K_C^{1/2} \rightarrow \big(K_C^{-1/2}\big)^*$ and $K_C^{-1/2} \rightarrow \big(K_C^{1/2}\big)^*$.

\section{The hyperk\"ahler structure of the Hitchin moduli space}\label{sec:2}

As a hyperk\"ahler space, the Hitchin moduli space has a rich geometric structure. At least two ongoing lines of research motivate us to consider the hyperk\"ahler geometry.
\begin{itemize}\itemsep=0pt
 \item There are many recent results about ``branes'' in Hitchin moduli space.
 \item There are recent results about the asymptotic geometry of the hyperk\"ahler metric.
\end{itemize}
In this section, we give an introduction to hyperk\"ahler geometry before specializing to the hyperk\"ahler geometry of the Hitchin moduli space. An excellent additional reference is~\cite{Neitzkeclass}.

\subsection{Introduction to hyperk\"ahler geometry}

A hyperk\"ahler manifold is a manifold whose tangent space admits an action of $I$, $J$, $K$ compatible with a single metric. To give a more precise definition of ``hyperk\"ahler'', we first have to define ``K\"ahler''.

\begin{defn} Let $(X,I)$ be a complex manifold of $\dim_\C X=n$.
 \begin{itemize}\itemsep=0pt
 \item A \emph{hermitian metric} on $(X, I)$ is a Riemannian metric $g$ such that $g(v,w) = g(Iv, Iw)$.
 \end{itemize}
 Let $\nabla$ denote the Levi-Civita connection\footnote{Recall that the Levi-Civita connection on $(X,g)$ is the unique connection that (1) preserves the metric, i.e., $\nabla g =0$ and (2) is torsion-free, i.e., for any vector fields $\nabla_X Y - \nabla_Y X =[X,Y]$.} on $TX$ induced by~$g$.
 \begin{itemize}\itemsep=0pt
 \item A hermitian metric $g$ on $(X, I)$ is \emph{K\"ahler} if $\nabla I=0$.
 \item If $(X, g, I)$ is K\"ahler, the \emph{K\"ahler form} $\omega \in \Omega^{1,1}(X)$ is defined by $\omega(v,w) = g(Iv, w)$.
 \end{itemize}
\end{defn}

\begin{ex} $\C$ is K\"ahler. The complex structure $I$ is given by multiplication by ${\rm i}$, the Riemannian metric is
 $g=\de x^2 + \de y^2 = \de z \de \zbar$,
 and the K\"ahler form is $\omega = \de x \wedge \de y = \frac{{\rm i}}{2} \de z \wedge \de \zbar$.
\end{ex}
It is easy to find examples of K\"ahler manifolds. For example, any complex submanifold of $\CP^n$ inherits a K\"ahler metric. Hyperk\"ahler manifolds are much more rigid, so it is harder to find examples.

\begin{defn} A \emph{hyperk\"ahler manifold} is a tuple $(X, g, I,J,K)$ where $(X,g)$ is a Riemannian manifold equipped with $3$ complex structures $I$, $J$, $K$~-- obeying the usual quaternionic relations~-- such that
 $(X, g, \bullet)$ is K\"ahler, for $\bullet=I,J,K$.
\end{defn}
The complex structures $I$, $J$, $K$ fit together into a $\zeta \in \CP^1$-family of complex structures
\begin{gather*}
 I_\zeta = \frac{1-|\zeta|^2}{1+|\zeta|^2} I + \frac{\zeta + \overline{\zeta}}{1+|\zeta|^2} J - \frac{{\rm i}(\zeta-\overline{\zeta})}{1+|\zeta|^2} K.
\end{gather*}
Consequently, given a hyperk\"ahler manifold $X$, we have a $\zeta \in \CP^1$-family of K\"ahler manifolds $(X, g, I_\zeta, \omega_\zeta)$.

\begin{ex}The vector space of quaternions, $\mathbb{H}$, is hyperk\"ahler
 \begin{align*}
 \mathbb{H} & \rightarrow \R^4,\\
 x_0 + x_1 {\rm i} + x_2 {\rm j} + x_3 {\rm k} & \mapsto (x_0, x_1, x_2, x_3).
 \end{align*}
 The three complex structure $I$, $J$, $K$ are respectively given by multiplication by ${\rm i}$, ${\rm j}$, ${\rm k}$. The hyperk\"ahler metric is $g= \de x_0^2 + \de x_1^2 + \de x_2^2 + \de x_3^2$. The symplectic forms $\omega_I$, $\omega_J$, $\omega_K$ are deter\-mi\-ned by the K\"ahler condition; e.g., since $\omega_I(v, w)=g(Iv,w)$, $\omega_I = \de x_0 \wedge \de x_1 + \de x_2 \wedge \de x_3$.
\end{ex}

\begin{exercise} \label{exercise:symplecticforms}Compute the symplectic forms $\omega_J, \omega_K \in \Omega^2\big(\R^4\big)$ for quaternion space.
\end{exercise}

As further evidence of the relevance of the quaternions for hyperk\"ahler manifolds, note that all hyperk\"ahler manifolds are automatically of dimension $\dim_\R X = 4k$ for $k \in \mathbb{N}$; furthermore,
$g$ is hyperk\"ahler if, and only if, its holonomy $\operatorname{Hol}_\nabla$ is a subgroup of $\operatorname{Sp}(k)$, the group of $k \times k$ quaternionic unitary matrices.

\subsection[Classification of noncompact hyperk\"ahler manifolds $X^4$]{Classification of noncompact hyperk\"ahler manifolds $\boldsymbol{X^4}$} \label{sec:type}
The hyperk\"ahler metric on the Hitchin moduli space is expected to be of type ``quasi-ALG'' a~generalization of ``ALG''. In this section, we explain the terminology ``ALG'' by taking a (somewhat lengthy) detour into the classic classification of noncompact $4$-dimensional hyperk\"ahler manifolds.

A noncompact complete connected hyperk\"ahler manifold $X$ of real dimension $4$ is called a \emph{gravitational instanton}\footnote{In some definitions of ``gravitational instanton'' the bound on the curvature is weakened. See Remark~\ref{rem:otherdefn}.} if there is some $\eps>0$ such that the Riemannian curvature tensor $\mathrm{Rm}$ satisfies the bound
\begin{gather}\label{eq:curvaturedecay}
 |\mathrm{Rm}|(x) \leq r(x)^{-2 -\eps},
\end{gather}
for $x \in X$ where $r(x)$ denotes the metric distance to a base point $o$ in $X$~\cite{ChenChen}.
Gravitational instantons can be divided into four categories:\\
\begin{tabular}{lll l}
$\bullet$ \textbf{ALE} & ``asymptotically locally Euclidean'' & $O\big(r^4\big)$ \qquad & ex) $\mathbb{H} \simeq \R^4$,\\
$\bullet$ \textbf{ALF} & ``asymptotically locally flat'' & $O\big(r^3\big)$ & ex) $ \R^3 \times S^1$,\\
$\bullet$ \textbf{ALG} & [\textsc{not an abbreviation}] & $O\big(r^2\big)$ & ex) $\R^2 \times T^2 $,\\
$\bullet$ \textbf{ALH} & [\textsc{not an abbreviation}] & $O\big(r^1\big)$ & ex) $\R \times T^3$.
\end{tabular}\\
Here, this coarse classification is by the dimension of the asymptotic tangent cone. The asymptotic tangent cones are,
respectively,\\
\begin{tabular}{l l}
$\bullet$ \textbf{ALE} & $\C^2/\Gamma$ where $\Gamma$ is a finite subgroup of $\mathrm{SU}(2)$,\\
$\bullet$ \textbf{ALF} & $\R^3$ or $\R^3/\Z_2$,\\
$\bullet$ \textbf{ALG} & $\C_\beta$, where $\C_\beta$ is a cone of angle $2\pi \beta$ for $\beta \in (0,1]$,\\
$\bullet$ \textbf{ALH} & $\R^+$.
\end{tabular}\\
Within each broad category (ALE/ALF/ALG/ALH), we have a finer classification by geometric type. Chen--Chen proved that any connected complete gravitational instanton with curvature decay like~\eqref{eq:curvaturedecay} must be asymptotic to some standard model. This is the data of a geometric type~\cite{ChenChen}. For each geometric type, we have a moduli space of hyperk\"ahler manifolds of that type. We will focus on the geometric classification for the cases ALE and ALG, since the ALE story is classical and the ALG story is most relevant for the Hitchin moduli space.

The geometric classification of ALE hyperk\"ahler metrics has been completed. The data for the geometric type is a finite subgroup~$\Gamma$ of~$\mathrm{SU}(2)$. Using this subgroup, define the singular space
\begin{gather*}X_\Gamma^\circ = \mathbb{C}^2/\Gamma.
\end{gather*}
Every ALE hyperk\"ahler $4$-manifold is diffeomorphic to the minimal resolution of~$X_\Gamma^\circ$ for some~$\Gamma$ \cite{KronheimerTorelli}. The moduli space $\cM_\Gamma$ of ALE instantons of type $\Gamma$ is non-empty and is parameterized by the integrals of the K\"ahler forms $\omega_I$, $\omega_J$, $\omega_K$ over the integer-valued second-homology lattice~\cite{Kronheimerconstruction, KronheimerTorelli}. The asymptotic tangent cone of any $X \in \cM_{\Gamma}$ is $\C^2/\Gamma$.

\begin{ex}$\Gamma= \Z_k$ acts on $(z=x_0 + {\rm i} x_1, w =x_2 + {\rm i} x_3)$ by $(z, w ) \mapsto \big(\e^{2 \pi \I/k} z, \e^{2 \pi \I/k} w\big)$. The moduli space $\cM_{\Z_k}$ has dimension $3k-6$ \cite{GibbonsHawking, HitchinALE}.
\end{ex}

For ALG gravitational instantons, the finer geometric classification is by the geometry at infinity. These standard models are torus bundles over the flat cone $\C_\beta$ of cone angle $2\pi \beta \in (0, 2\pi]$. The list of torus bundles $E \rightarrow \C_\beta$ is quite restricted.
\begin{defnthm}[{\cite[Theorem 3.11]{ChenChen}, \cite[Theorem 3.2]{ChenChenIII}}] Suppose $\beta \in (0, 1]$ and $\tau \in \mathbb{H}=\{ \tau \,|\, \mathrm{Im}(\tau) >0\}$ are parameters in the following table:
\begin{gather}\label{eq:table}
\begin{array}{| l | c | c | c | c | c | c | c | c |}
\hline
D & \mbox{Regular} & I_0^* & II & II^* & III & III^* & IV & IV^* \\
\beta & 1 & \frac{1}{2} & \frac{1}{6} & \frac{5}{6} & \frac{1}{4} & \frac{3}{4} & \frac{1}{3} & \frac{2}{3} \\
\tau & \in \mathbb{H} & \in \mathbb{H} & \e^{2 \pi \I/3} & \e^{2 \pi \I/3} & \I & \I & \e^{2 \pi \I/3} & \e^{2 \pi \I/3}
\\\hline
\end{array}
\end{gather}
Suppose $\ell>0$ is some scaling parameter. Let $E$ be the manifold obtained by identifying the two boundaries of the torus bundle over the sector
\begin{gather*}
\{u \in \C\colon \mathrm{Arg}(u) \in [0, 2\pi \beta] \; \& \; |u| \geq R\} \times \C_v/(\Z \ell + \Z \ell \tau)
\end{gather*}
by the gluing map $(|u|,v) \simeq \big(\e^{2 \pi \I \beta}|u|, \e^{2 \pi \I \beta}v\big)$. This manifold together with a certain (see \cite[Definition~2.3]{ChenChenIII}) flat hyperk\"ahler metric $g_{\mathrm{mod}}$ is called the \emph{standard ALG model of type $(\beta, \tau)$}.

Every ALG gravitational instanton $X$ is asymptotic to the one of these standard models $(E, g_{\mathrm{model}})$. Moreover, if $\beta = 1$, then $X$ \emph{is} the standard flat gravitational instanton $\C \times T^2_\tau$.
\end{defnthm}

To explain why Kodaira types of singular fibers appear in the first row of \eqref{eq:table}, note the following theorem:
\begin{thm}[\cite{ChenChen}] \label{thm:ALGfibration} Any ALG gravitational instanton $X$ can be compactified in a complex analytic sense. I.e., there exists a compact elliptic surface\footnote{Furthermore, from \cite[Theorem~1.3]{ChenChenIII}, $\psi\colon \overline{X} \rightarrow \CP^1$ is a rational elliptic surface in the sense of \cite[Definition~2.7]{ChenChenIII}} $\overline{X}$ with a meromorphic function
$\psi\colon \overline{X} \rightarrow \CP^1$ whose generic fiber is a complex torus. The fiber $D=\psi^{-1}(\infty)$ is either regular or singular of Kodaira type $I^*_0$, $II$, $II^*$, $III$, $III^*$, $IV$, $IV^*$. Moreover, there is some $\zeta \in \CP^1$ such that $(X, I_\zeta)$ is biholomorphic to $\overline{X}-D$.
\end{thm}

\begin{rem} Looking forward, this fibration $\psi$ should loosely remind you of the Hitchin fibration $\operatorname{Hit}\colon \cM \rightarrow \cB \simeq \C^{\frac{1}{2} \dim_\C \cM}$.
\end{rem}

\begin{rem} \label{rem:otherdefn}There are other definitions of gravitational instantons appearing in the literature without the strict curvature bounds in \eqref{eq:curvaturedecay}. Stranger things can happen if we remove these curvature bounds. Given any noncompact complete connected hyperk\"ahler manifold $X$ of real dimension $4$, one can associate a number based on the asymptotic volume growth of $B_r$, a ball of radius $r$. With the definition of gravitational instantons in \eqref{eq:curvaturedecay}, the volume growth is an integer: 4, 3, 2, 1. Without the curvature bounds in~\eqref{eq:curvaturedecay}, the volume growth need not be an integer. There are no hyperk\"ahler metrics with growth between $r^3$ and $r^4$ \cite{minerbe}. However, Hein constructed an example of a hyperk\"ahler metric with volume growth~$r^{4/3}$~\cite{hein43}. Chen--Chen call this an example of type ALG$^*$ since the growth rate $\frac{4}{3}$ is the in ALG-like interval $(1, 2]$. The~``$*$'' indicates that the modulus of the torus fibers is changing; it is unbounded asymptotically, i.e., the torus fiber is becoming very long and thin.
\end{rem}

\subsection{The hyperk\"ahler metric on the Hitchin moduli space}\label{sec:metric}

The Hitchin moduli space has a hyperk\"ahler metric.
To hint at the origins of the hyperk\"ahler structure,
we instead discuss the origin of the $\CP^1$-family of complex structures on the Hitchin moduli space. The $\CP^1$-family of complex structures arises from the
 complex structure on the Riemann surface $C$ and the complex structure on the group $G_\C$; these respectively, give the $I$ and $J$ complex structures on $\cM$.

 Each of these complex structures $I_\zeta$ gives an avatar of the Hitchin moduli space
 as a complex manifold:
\begin{itemize}\itemsep=0pt
 \item $\cM_{\zeta =0}=(\cM, I_{\zeta=0})$ is the Higgs bundle moduli space;
 \item $\cM_{\zeta \in \C^\times}$ is the moduli space of flat connections;
 \item $\cM_{\zeta=\infty}$ is the moduli space of anti-Higgs bundles.
\end{itemize}
Note that these can be genuinely different as complex manifolds. For $\cM=\cM\big(T^2_\tau, \GL(1,\C)\big)$,
\begin{gather} \label{eq:ex}
 \cM_0 \simeq \C \times T^2_\tau, \qquad \cM_{\zeta \in \C^\times} \simeq \C^\times \times \C^\times, \qquad \cM_\infty \simeq \C \times T^2_{- \overline{\tau}},
\end{gather}
where $T^2_\tau = \C/(\Z \oplus \tau \Z)$ is the complex torus with parameter $\tau$. The hyperk\"ahler metric on $\cM_0 \simeq \R^2_{x_0, x_1} \times T^2_{x_3, x_4}$ is $g=\de x_0^2 + \de x_1^2 + \de x_2^2 + \de x_3^2$, which is indeed $ALG$.

In general, the hyperk\"ahler metric on Hitchin moduli space is expected to be of type ``quasi-ALG'' which is some generalization\footnote{This term ``QALG'' has not been formally defined, but QALG is supposed to generalize ALG in an analogous way as QALE generalizes ALE and QAC generalizes AC.

Dominic Joyce considered a higher-dimensional version of ALE, and subsequently defined QALE. In this context, a (Q)ALE metric is a K\"ahler metric on a manifold of real-dimension $2n$ with asymptotic volume growth like $r^{2n}$. Fix a finite subgroup $\Gamma \subset \mathrm{U}(n)$. If $\Gamma$ acts freely on $\C^n-\{0\}$, then $\C^n/\Gamma$ has an isolated quotient singularity at $0$. The appropriate class of K\"ahler metrics on the resolution $X$ of $\C^n/\Gamma$ \cite{JoyceALE} are ALE metrics. If however, $\Gamma$~does not act freely, then the singularities of $\C^n/\Gamma$ extend to the ends. The appropriate class of K\"ahler metrics on resolution $X$ of non-isolated quotient singularities are called quasi-ALE or QALE~\cite{JoyceQALE}.

For quasi-asymptotically conical (QAC) versus AC, see for example~\cite{QAC}.\label{rem:QALG}} of ALG. In higher dimensions ``ALG'' has not be formally defined; however, by analogy with the $4$-dimensional case described in Section~\ref{sec:type}, any higher-dimensional $ALG$ hyperk\"ahler manifold should be asymptotic to a flat torus bundle over a~half-dimensional complex vector space. Moreover, the modulus of the torus lattice should stay bounded. (This condition on the modulus rules out higher-dimensional generalizations of Hein's ALG$^*$ example in Remark~\ref{rem:otherdefn}.)

For the Hitchin moduli space, the Hitchin fibration $\operatorname{Hit}\colon \cM \rightarrow \cB$ should asymptotically give the torus fibration over a half-dimensional complex vector space. It certainly does in the ALG example of~\eqref{eq:ex}! Because the singular locus typically intersects the ends of the Hitchin base $\cB$, we will typically not have a nondegenerate asymptotic torus fibration $\psi\colon \cM \rightarrow \cB$. Consequently, in these cases, the Hitchin moduli space is instead expected to be ``QALG'', as described in footnote~\ref{rem:QALG}, rather than ALG.

\subsection{The hyperk\"ahler metric} \label{sec:hk}

Alternatively, the Hitchin moduli space can be viewed as pairs $\big[(\delbar_A, \Phi)\big]$ solving Hitchin's equations given a fixed complex vector bundle $E \rightarrow C$ \emph{with} fixed hermitian metric. Once we've fixed a hermitian metric, we only consider gauge transformations which fix the hermitian metric. This gives a reduction from complex gauge transformations $\mathcal{G}_\C=\mathcal{GL}(E)$ to unitary gauge transformations $\mathcal{G}=\mathcal{U}(E)$-gauge transformations.

The hyperk\"ahler metric on $\cM$ is defined using the unitary formulation of Hitchin's equations in terms of pairs $(\delbar_A, \Phi)$. First, consider the configuration space $\mathcal{C}$ of all pairs $(\delbar_A, \Phi)$ solving Hitchin's equations~-- without taking gauge equivalence. The space of holomorphic structures on $E$ is an affine space modeled on $\Omega^{0,1}(C, \End E)$. The space of all Higgs fields is the vector space $\Omega^{1,0}(C, \End E)$. Thus, the set of pairs $(\delbar_A, \Phi)$ solving Hitchin's equations sits inside an affine space modeled on
\begin{gather*}
 \Omega^{0,1}(C, \End E) \times \Omega^{1,0}(C, \End E).
\end{gather*}
This product space has a natural $L^2$-metric given by
\begin{gather*}
 g\big(\big(\dot{A}^{0,1}_1, \dot\Phi_1\big), \big(\dot{A}^{0,1}_2, \dot\Phi_2\big) \big) =2\I \int_C \big\langle \dot{\Phi}_1 \,\overset{\wedge}{,}\, \dot{\Phi}_2\big\rangle - \big\langle \dot{A}^{0,1}_1 \,\overset{\wedge}{,} \,\dot{A}^{0,1}_2\big\rangle,
\end{gather*}
where the hermitian inner products are taken only on the matrix-valued piece, so that $\IP{\dot{\Phi}_1 \,\overset{\wedge}{,}\, \dot{\Phi}_2}$ is a $(1,1)$-form. The factor $2\I$ appears since $\de z \wedge \de \zbar = - 2\I \de x \wedge \de y$. The hyperk\"ahler metric $g_{\cM}$ on $\cM$ descends from this $L^2$-metric $g$. Note that any tangent vector $[(\dot{A}^{0,1}, \dot{\Phi})] \in T_{[(\delbar_A, \Phi)]} \cM$ has multiple representatives. The hyperk\"ahler metric $g_{\cM}$ on $\cM$ is defined so that
\begin{gather*}
 \norm{[(\dot{A}^{0,1}, \dot{\Phi})]}_{g_{\cM}} = \min_{(\dot{A}^{0,1}, \dot{\Phi}) \in [(\dot{A}^{0,1}, \dot{\Phi})]} \quad \norm{(\dot{A}^{0,1}, \dot{\Phi}) }_g.
\end{gather*}
The minimizing representative is said to be in ``Coulomb gauge''. We will call this natural hyperk\"ahler metric $g_\cM$ ``Hitchin's hyperk\"ahler $L^2$-metric''.

\subsection{Branes in the Hitchin moduli space}
Recently, there have been a number of results about branes in the Hitchin moduli space. (See~\cite{3Lauras} for a survey of results and further directions.)

\begin{defn} A \emph{brane} is an object in one of the following categories:
 \begin{itemize}\itemsep=0pt
 \item[$\bullet$] [$A$-side, i.e., symplectic] Fukaya category, or
 \item[$\bullet$] [$B$-side, i.e., complex] derived category of coherent sheaves.
 \end{itemize}
 \end{defn}
\noindent The approximate data of an $(A/B)$-brane in a (symplectic/complex) manifold $X$ is
\begin{itemize}\itemsep=0pt
 \item a submanifold $Y \subset X$, together with
 \item $(E, \nabla) \rightarrow Y$, a vector bundle with connection.
\end{itemize}
Further ignoring bundles, an \emph{A-brane} in a symplectic manifold $(X, \omega)$ ``is'' a Lagrangian submanifold\footnote{Let $(X, \omega)$ be a symplectic manifold. Recall a submanifold $L$ is Lagrangian if $\dim_\R L = \frac{1}{2} \dim_R X$ and $\omega|_L =0$.} of $X$. A \emph{B-brane} in a complex manifold $(X,I)$ ``is'' a holomorphic submanifold.

If $X$ is hyperk\"ahler, then $X$ has a triple of K\"ahler structures $(X,g, I, \omega_I, J, \omega_J, K, \omega_K)$. With respect to the triple of K\"ahler structures, a \emph{$(B, A, A)$ brane} in $X$ ``is'' a submanifold $Y$ which is holomorphic with respect to $I$, Lagrangian with respect to $\omega_J$, and Lagrangian with respect to~$\omega_K$. Not all types exis~-- only $(B,A,A), (A,B,A), (A,A,B)$ and $(B,B,B)$-branes exist.

\begin{figure}[ht]\centering
\includegraphics[height=1.0in]{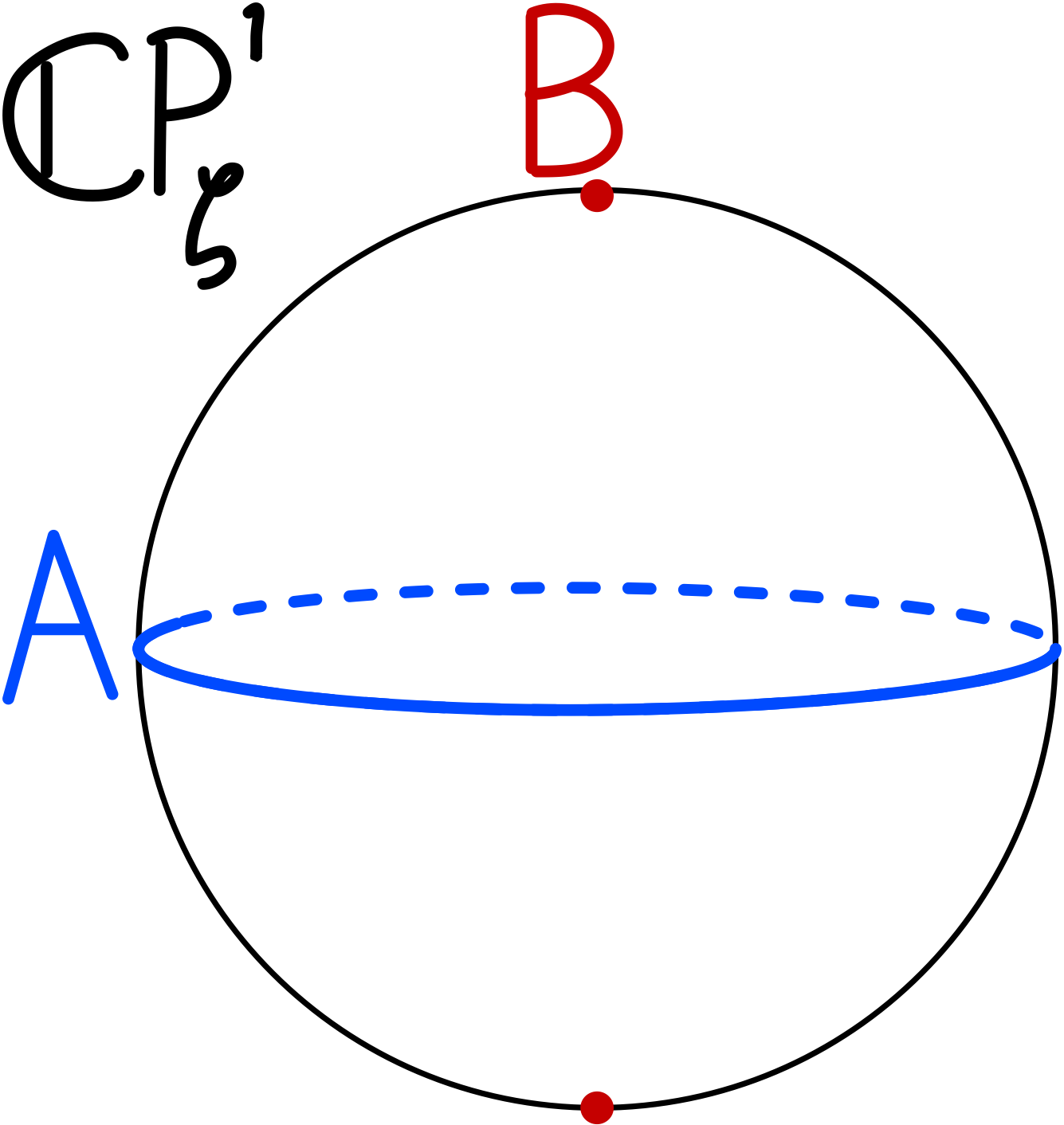}
\caption{$(B,A,A)$-brane.}\label{fig:branes}
\end{figure}

Given a $(B, A, A)$ brane $Y\subset X$, one might ask whether the submanifold $Y$ is holomorphic with respect to $I_\zeta$ or Lagrangian with respect to $\omega_\zeta$ for any of the other K\"ahler structures $(X, g, I_\zeta, \omega_\zeta)$. In fact, as shown in Fig.~\ref{fig:branes}, $Y$ is holomorphic with respect to both $\pm I$, and~$Y$ is Lagrangian for $|\zeta|=1$~-- the whole circle containing~$J$ and~$K$. It is neither holomorphic, nor Lagrangian for any other value $\zeta$. Similar statements hold for $(A,B,A)$ and $(A,A,B)$ branes.

\begin{exercise}Let $Y$ be the $(x_0, x_1)$-plane in quaternion space.
\begin{itemize}\itemsep=0pt
\item[(a)]Show that $Y$ is (the support of) a $(B,A,A)$ brane, i.e., it's holomorphic with respect to $I$ and Lagrangian with respect to $\omega_J$ and $\omega_K$.
\item[(b)] Is $Y$ holomorphic with respect to any other complex structure $I_\zeta$? Is $Y$ Lagrangian with respect to any other symplectic structure $\omega_\zeta$?
\end{itemize}
\end{exercise}

Many recent results concern constructions of different families of branes inside the Hitchin moduli space. For example, some branes~-- including the Hitchin section, which is a $(B,A,A)$-brane~-- are constructed as fixed point sets of certain involutions on the moduli space of Higgs bundles. (See, for example,~\cite{slices}). Langlands duality, shown in~\eqref{eq:Langlands}, exchanges the brane types. For example, $(B, A, A)$-branes in the $G$-Hitchin moduli space $\cM_{G}$ get mapped to $(B,B,B)$-branes in the $^L G$-Hitchin moduli space~$\cM_{^L G}$
\begin{equation} \label{eq:Langlands}
\begin{tikzcd}
\cM_{G}
\arrow{rd} &
&\cM_{^L G} \arrow{ld}\\
& \cB_{G} \simeq \cB_{^L G}. &
\end{tikzcd}
\end{equation}

\section{Spectral interpretation and limiting configurations}\label{sec:3}

In Section~\ref{sec:spectraldata}, we introduced the Hitchin fibration. Now, we
\begin{enumerate}\itemsep=0pt
 \item[1)] give a geometric interpretation of the Hitchin fibration, and
 \item[2)] give a construction of the harmonic metric for a Higgs bundle near the ends
 of the moduli space.
\end{enumerate}

\subsection{Spectral data}
The Hitchin fibration, introduced in \eqref{eq:Hitchinfibration}, is a map
\begin{align*}
 \operatorname{Hit}\colon \cM &\twoheadrightarrow \cB \simeq \C^{\frac{1}{2} \dim_\C \cM},\\
 \big(\delbar_E, \varphi, h\big) &\mapsto \operatorname{char}_\varphi (\lambda),
\end{align*}
where the characteristic polynomial $\operatorname{char}_\varphi\lambda$
encodes the eigenvalues $\lambda_1, \dots, \lambda_n$ of $\varphi$. There are two additional interpretations of $\cB$ that are useful:
\begin{figure}[ht]\centering
\includegraphics[height=1.0in]{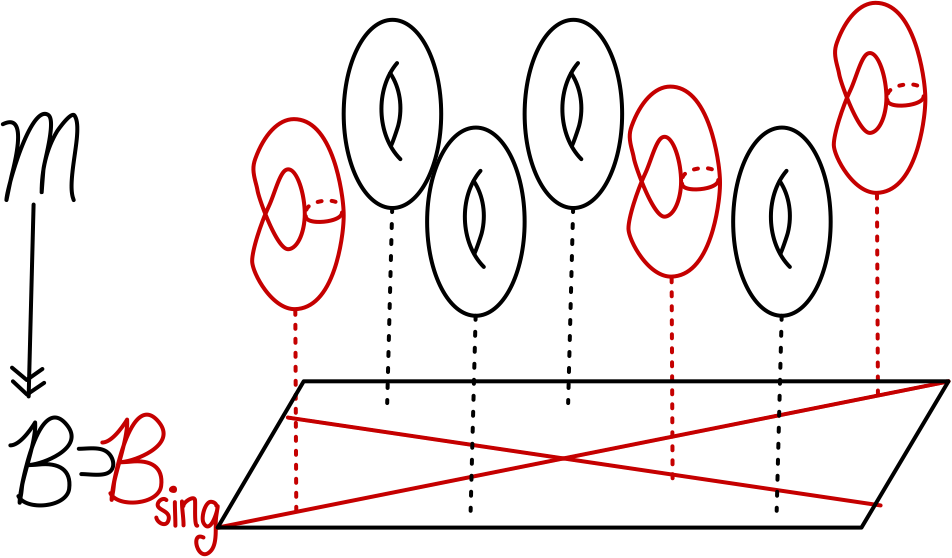}
\caption{Hitchin fibration.}\label{fig:hitchinfibration}
\end{figure}
\begin{itemize}\itemsep=0pt
 \item [$\bullet$] \textbf{[algebraic interpretation]} The coefficients of
 \begin{gather*}\operatorname{char}_\varphi (\lambda)= \lambda^n + q_1 \lambda^{n-1} + q_2 \lambda^{n-2} + \cdots + q_{n-1} \lambda + q_n
 \end{gather*}
 are sections $q_i \in H^0\big(C, K_C^i\big)$. Consequently, the Hitchin base $\cB$ can be identified with the complex vector space of coefficients of $\operatorname{char}_\varphi (\lambda)$. For example,
\begin{gather*}
 \cB_{\GL(n,\C)} = \bigoplus_{i=1}^n H^0\big(C, K_C^i\big) \ni (q_1, \ldots, q_n), \qquad \mbox{and}\\
 \cB_{{\rm SL}(n,\C)} = \bigoplus_{i=2}^n H^0\big(C, K_C^i\big) \ni (q_2, \ldots, q_n),
\end{gather*}
since for ${\rm SL}(n,\C)$, $q_1 = -\operatorname{tr} \varphi =0$.
\item [$\bullet$]\textbf{[geometric interpretation]}
Define
\begin{gather*}
 \Sigma =\{ \lambda\colon \operatorname{char}_\varphi (\lambda) =0\} \subset \operatorname{Tot}(K_C).
\end{gather*}
Then $\Sigma$ is ``spectral cover''.
The spectral cover $\pi\colon \Sigma \rightarrow C$, shown in Fig.~\ref{fig:spectralcover}, is a ramified $n:1$ cover of $C$. The branch locus $Z$ is the zero locus of the discriminant section
\begin{gather*} %\label{eq:discriminant}
 \Delta_\varphi = \prod_{i<j} (\lambda_i - \lambda_j)^2.
\end{gather*}
\begin{figure}[ht]\centering
\includegraphics[height=1.0in]{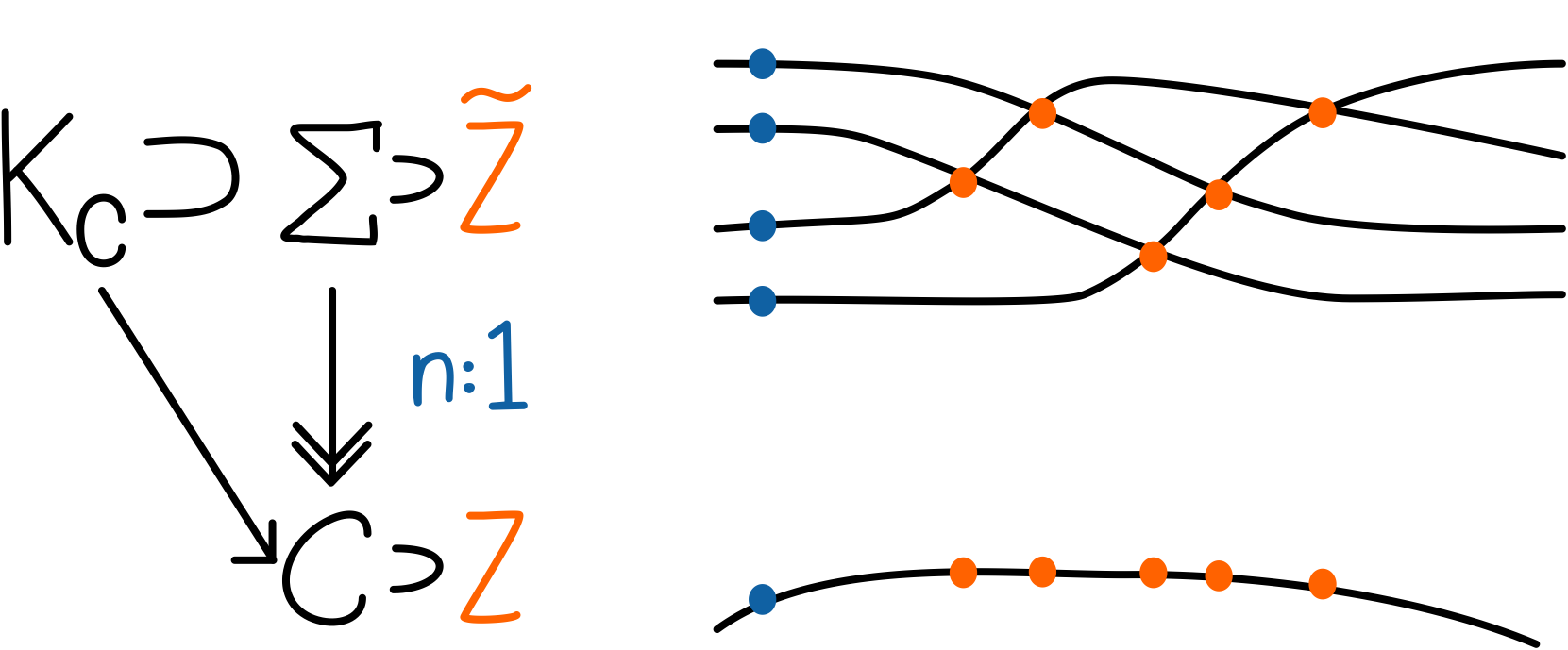}
\caption{Each of the (generically) $n$ sheets of $\Sigma$ represents an eigenvalue of $\varphi$. The branch locus $Z$ is shown in orange.}\label{fig:spectralcover}
\end{figure}
\end{itemize}

For $\GL(n,\C)$ and ${\rm SL}(n,\C)$-Higgs bundles, the regular locus $\cM'$, discussed in Section~\ref{sec:spectraldata}, consists of Higgs bundles lying over \emph{smooth} spectral covers $\Sigma$.

\begin{exercise} \quad
\begin{itemize}\itemsep=0pt
 \item[(a)] What bundle is the discriminant section $\Delta_\varphi$ a section of? What is the number of zeros of~$\Delta_\varphi$ with multiplicity?
\item[(b)]Use this to compute the number of ramification points of $\pi\colon \Sigma \rightarrow C$ (with multiplicity).
\item[(c)] Compute the genus of $\Sigma$.\\
 \texttt{Hint:} \emph{In the case of unramified $N:1$ covers $\pi\colon S' \rightarrow S$, the Riemann--Hurwitz formula says that $\chi(S')=N \chi(S)$ where $\chi(S)=2(\gamma_S-1)$ is the Euler characteristic and $\gamma_S$ is the genus of~$S$. In the case of ramified covers, this is corrected to}
 \begin{gather*}
 \chi(S')=N \chi(S) - \sum_{P \in \tilde{Z}} (e_P-1),
 \end{gather*}
 \emph{where $\tilde{Z} \subset S'$ is the set of points of $S'$ where $\pi$ locally looks like $\pi(z)=z^{e_P}$. The number~$e_P$ is called the ``ramification index''.}
\end{itemize}
\end{exercise}

\begin{fact}For ${\rm SL}(2,\C)$, $\operatorname{Hit}\colon \big(\delbar_E, \varphi, h\big) \mapsto \det \varphi =q_2$, so the branch locus $Z$ is the set of zeros of $\det \varphi$. Call $p \in Z$ a simple zero if $q_2 \sim z \de z^2$, and call $p$ a $k^{\rm th}$ order zero if $q_2 \sim z^k \de z^2$. For $\cM_{{\rm SL}(2,\C)}$, $q_2$ has only simple zeros $\Leftrightarrow$ $\Sigma$ is smooth $\Leftrightarrow$ the spectral cover $\Sigma$ lies in the regular locus $\cB' = \cB-\cB_{{\rm sing}}$ $\Leftrightarrow$ $\operatorname{Hit}^{-1}(\Sigma)$ is a compact abelian variety.
\end{fact}

\begin{figure}[ht]\centering
\includegraphics[height=1.0in]{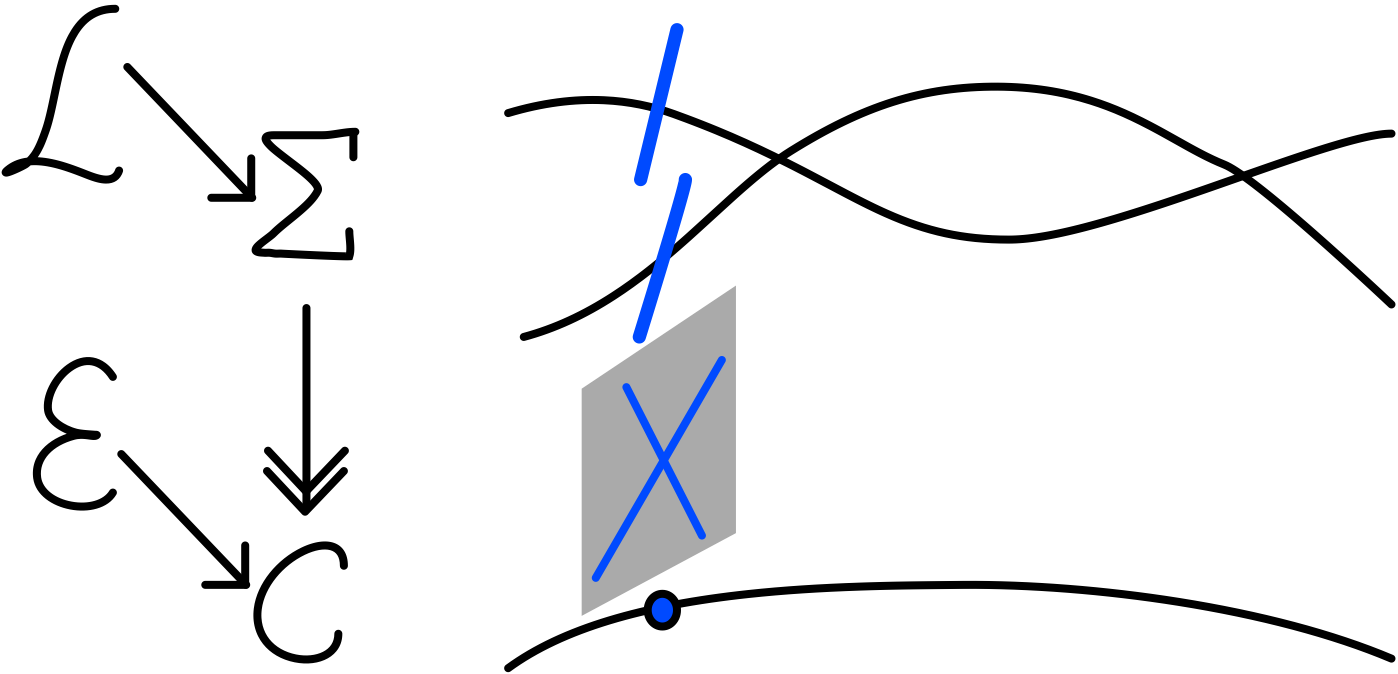}
\caption{Each sheet of the spectral cover $\Sigma \to C$ over a point $x\in C$ corresponds to an eigenvalue of~$\varphi(x)$. The fiber of the spectral line bundle $\cL$ is the associated eigenspace of $\varphi(x)$.}\label{fig:line}
\end{figure}

\looseness=1 Having given a geometric interpretation of the Hitchin base~$\cB$, we now give a geometric interpretation of the torus fibers of $\operatorname{Hit}\colon \cM' \to \cB'$. (We restrict to the regular locus \smash{$\cM' \rightarrow \cB'$} since the torus fibers degenerate over the singular locus $\cB_{{\rm sing}}$.)
As shown in Fig.~\ref{fig:line}, the eigenspaces of a Higgs field $\varphi$ can be encoded in a line bundle $\cL \rightarrow \Sigma$. Note that $\cE \simeq \pi_*\cL$. The torus fiber $\operatorname{Hit}^{-1}(\Sigma)$ is some space of line bundles over the spectral cover~$\Sigma$. For $\GL(n,\C)$, this fiber is the Jacobian, $\operatorname{Jac}(\Sigma)$. For ${\rm SL}(n,\C)$, this fiber is the Prym variety $\operatorname{Prym}(\Sigma, C)$. Here, we have a~subvariety of the Jacobian because of the trivialization of determinant as $\Det (\pi_* \cL) \simeq \Det \cE \simeq \cO$.

\subsection{Limits in the Hitchin moduli space}\label{sec:limits}

In \cite{GMNwallcrossing, GMNhitchin}, Gaiotto--Moore--Neitzke give a conjectural description of the hyperk\"ahler metric~$g_\cM$. Gaiotto--Moore--Neitzke's conjecture suggests that~-- surprisingly~-- much of the asymptotic geo\-met\-ry of~$\cM$ can be derived from the abelian data $\cL \rightarrow \Sigma$. (We give a survey of this conjecture and recent progress in Section~\ref{sec:4}.) In this section, we describe how solutions of Hitchin's equations \emph{at} the ends of the moduli space come naturally from the abelian data $\cL \rightarrow \Sigma$. This is an early hint of the importance of the abelian data for the asymptotic geometry of the Hitchin moduli space.

\begin{figure}[ht]\centering
\includegraphics[height=1.0in]{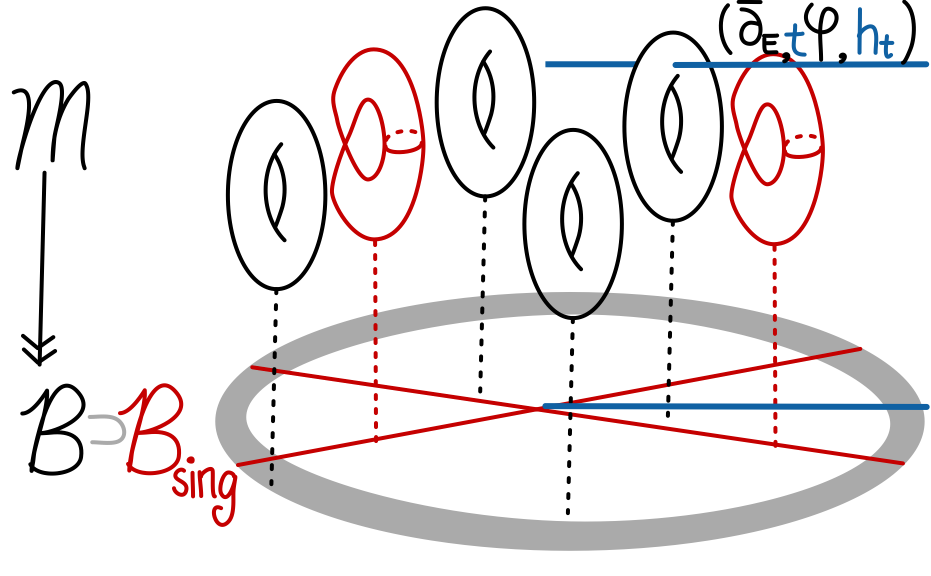}
\caption{A ray $\big(\delbar_E, t \varphi, h_t\big)$ going off to the $t=\infty$ ends of $\cM$.}\label{fig:ray}
\end{figure}

As shown in Fig.~\ref{fig:ray}, consider the ray $\big(\delbar_E, t \varphi, h_t\big)$ of solutions of Hitchin's equations in $\cM$.
As $t \rightarrow \infty$, the curvature $F_{D(\delbar_E, h_t)}$ (1) concentrates at the branch locus $Z$ and (2) vanishes everywhere else (exponentially in~$t$) \cite{FredricksonSLn, MSWW14, Mochizukiends, Taubes3}. Consequently, the limiting hermitian metric~$h_\infty$ is singular at~$Z$ and solves the decoupled Hitchin's equations
\begin{gather*}
 F_{D(\delbar_E, h_\infty)} = 0, \qquad [\varphi, \varphi^{*_{h_\infty}}]=0.
\end{gather*}
In fact, the limiting metric $h_\infty$ is actually a pushforward of a singular hermitian metric (with singularities at the ramification points)~$h_\cL$ on the spectral line bundle~$\cL \rightarrow \Sigma$. Consequently, we say that ``Hitchin's equations abelianize at the ends''.

To construct $h_\infty$, we start with a Higgs bundle $(\cE, \varphi) \in \cM'$. As shown in~\eqref{eq:steps}, there are 4 steps:
\begin{equation} \label{eq:steps}
\begin{tikzcd}
\cE \arrow{dd} & \cL \arrow{dd} & (\cL, \alpha_{\widetilde{p}}) \arrow{dd}\\
 \arrow{r}{\mbox{\textcircled{1}}} & \arrow{r}{\mbox{\textcircled{2}}} & \; \\
C & \Sigma & \Sigma\arrow{d}{\mbox{\textcircled{3}}}\\
h_\infty & & h_\cL \arrow{ll}{\mbox{\textcircled{4}}}
\end{tikzcd}
\end{equation}
\begin{enumerate}\itemsep=0pt
\item Let $\cL \rightarrow \Sigma$ be the associated spectral data, consisting of the line bundle $\cL$ and spectral cover $\Sigma$.
\item Equip $\cL$ with certain\footnote{In the case $G_\C={\rm SL}(2,\C)$, these weights are $-\frac{1}{2}$. For higher rank, see \cite{FredricksonSLn}.} parabolic weights $\alpha_{\widetilde{p}}$ at $\widetilde{p} \in \widetilde{Z}$. This makes $(\cL, \alpha_{\widetilde{p}})$ a parabolic line bundle over $\Sigma$.
\item Let $h_\cL$ be the Hermitian--Einstein metric on $\cL \rightarrow \Sigma$
which is ``adapted'' to the parabolic structure. The hermitian metric
$h_\cL$ solves the Hermitian--Einstein equation $F_{D(\delbar_\cL, h_\ell)} =0$.
Because $h_\cL$ is ``adapted'' to the parabolic structure, $h_\cL$ has a singularity at $\widetilde{p} \in \widetilde{Z}$ like $h_\cL \simeq |w|^{2 \alpha_{\widetilde{p}}}$.
\item Finally, $h_\infty$ is the orthogonal pushforward of $h_\cL$. The eigenspaces of $\varphi$ are mutually orthogonal with respect to $h_\infty$, and $h_\infty$ agrees with $h_\cL$ in each eigenspace of $\varphi$.
\end{enumerate}

\section[Some recent results about the asymptotic geometry of $\cM$]{Some recent results about the asymptotic geometry of $\boldsymbol{\cM}$}\label{sec:4}

In \cite{GMNwallcrossing, GMNhitchin}, Gaiotto--Moore--Neitzke conjecture that Hitchin's hyperk\"ahler metric solves an integral relation\footnote{The integral equation appears in \cite[equation~(4.8)]{Neitzkehyperkahler}
which is a survey of~\cite{GMNhitchin} aimed at mathematical audiences.}. As a consequence of their conjecture, Hitchin's hyperk\"ahler metric~$g_{\cM}$ admits an expansion in terms of a simpler hyperk\"ahler metric $g_{\semif}$, known as the ``semiflat metric''. The semiflat metric $g_{\semif}$ is smooth on $\cM'$ and exists because $\cM$ is a algebraic completely integrable system~\cite[Theorem~3.8]{Freedsemiflat}.
 Thus it is deeply related to the spectral data.

\begin{conj}[{weak form of Gaiotto--Moore--Neitzke's conjecture for $\cM_{\SU(2)}$ \cite[equation~(1.1)]{DumasNeitzke}}] Fix a Higgs bundle $(\delbar_E, \varphi)$ in $\cM'$ and let $q_2=\det \varphi$. Hitchin's hyperk\"ahler $L^2$-metric on~$\cM'$ admits an expansion as
\begin{gather} \label{eq:weakGMN}
 g_{\cM} = g_{\semif} + O \big(\e^{-4Mt} \big),
\end{gather}
where $M$ is the length of the shortest geodesic on the associated spectral cover $\Sigma$, measured in the singular flat metric~$\pi^*|q_2|$, normalized so that~$C$ has unit volume.
\end{conj}
\begin{figure}[ht]\centering
\includegraphics[height=1.0in]{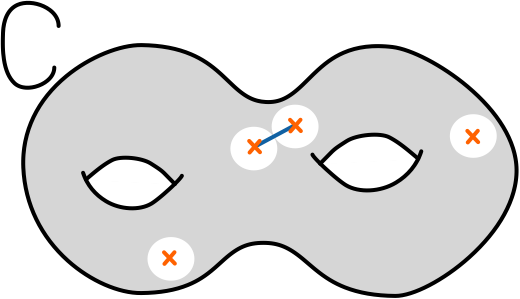}
\caption{The shortest geodesic is the lift of the indicated saddle connection between the zeros ``$\times$'' of~$q_2$.}\label{fig:min}
\end{figure}
\begin{rem}When we write $g_\cM = g_{\semif} + O\big(\e^{-4Mt}\big)$ in \eqref{eq:weakGMN} what we really mean is this: Fix a~Higgs bundle $(\cE, \varphi)$ in $\cM'$, and a Higgs bundle variation $\dot{\psi}=(\dot{E}, \dot{\varphi})$. Consider the deformation $\dot{\psi}_t=(\dot{E}, t \dot{\varphi}) \in T_{(\cE, t \varphi)} \cM$ over the ray $(\cE, t \varphi, h_t)$. As $t \rightarrow \infty$, the difference between
 Hitchin's hyperk\"ahler $L^2$-metric $g_{\cM}$ on $\cM$ and the semiflat (hyperk\"ahler) metric $g_{\semif}$ satisfies
 \begin{gather*} %\label{eq:summary}
 g_{\cM}(\dot{\psi}_t, \dot{\psi}_t) = g_{\semif}(\dot{\psi}_t, \dot{\psi}_t) + \cO\big(\e^{-4Mt}\big).
 \end{gather*}
\end{rem}

We briefly review the recent progress towards proving this conjecture for the $\SU(2)$-Hitchin moduli space, in chronological order. Mazzeo--Swoboda--Weiss--Witt \cite{MSWW17} have shown that along a generic ray $g_{\cM}-g_{\semif}$ decays \emph{polynomially} in~$t$. Dumas--Neitzke \cite{DumasNeitzke} have shown that~-- restricted to the Hitchin section~-- $g_{\cM}-g_{\semif}$ decays \emph{exponentially} in $t$ like $\cO\big(\e^{-2Mt}\big)$. The author has shown that along a ray in $\cM'_{\SU(2)}$ $g_{\cM}-g_{\semif}$ decays \emph{exponentially} in~$t$~\cite{Fredricksonasymptoticgeometry}. However, the constant of exponential decay is not sharp.

We now remark on a few ideas behind the proofs.

\textbf{Idea 1.} One crucial observation is the analogy
\begin{eqnarray*}
 \mbox{Hitchin's hyperk\"ahler $L^2$-metric $g_{\cM}$} & : & \mbox{the harmonic metric $h_t$} \\&::&\\
 \mbox{the semiflat metric $g_{\semif}$} &\; : \;& \mbox{the limiting metric $h_\infty$.}
\end{eqnarray*}
The regular locus of the Hitchin moduli space consists of triples $[(\delbar_E, \varphi, h)]$ where $h$ is the harmonic metric. Define the moduli space of limiting configurations to be the space of triples $[(\delbar_E, \varphi, h_\infty)]$; we replace the harmonic metric with the limiting metric $h_\infty$ from Section~\ref{sec:limits}. Recall that in Section~\ref{sec:metric}, that Hitchin's hyperk\"ahler $L^2$-metric $g_{\cM}$ was constructed as the $L^2$-metric on $\cM$. Similarly,

\begin{prop}[{\cite[Proposition 3.7, Proposition 3.11, Lemma 3.12]{MSWW17}}] The semiflat met\-ric~$g_{\semif}$ is the natural hyperk\"ahler $L^2$-metric on the moduli space of limiting configurations~$\cM'_\infty$, for deformations in Coulomb gauge.
\end{prop}

\begin{figure}[h!]\centering
\includegraphics[height=1.2in]{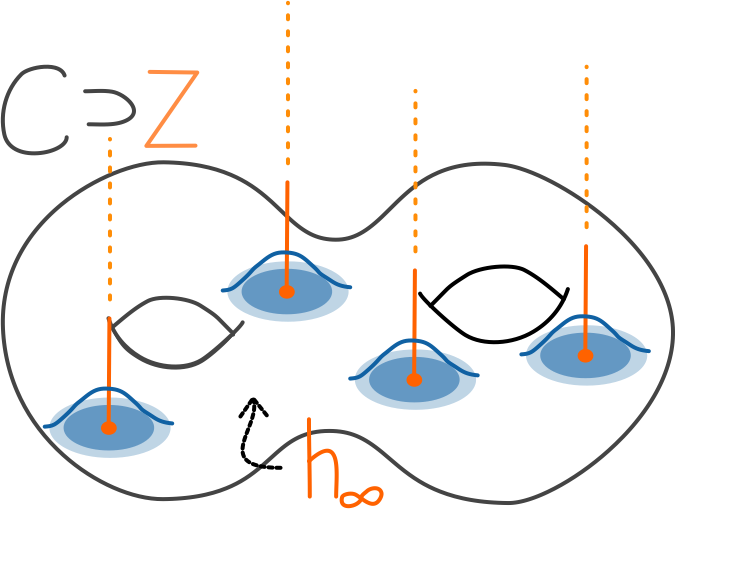}
\caption{Approximate solutions $h_t^\app$ are constructed by desingularizing $h_\infty$.}\label{fig:globalgluing}
\end{figure}

\textbf{Idea 2.} The result in \cite{MSWW17} is built on Mazzeo--Swoboda--Weiss--Witt's description of the harmonic metrics near the ends of the Hitchin moduli space in~\cite{MSWW14}. Mazzeo--Swoboda--Weiss--Witt build an family of approximate solutions of Hitchin's equations $\big(\delbar_E, t\varphi, h_t^\app\big)$ that are exponentially close to the actual solutions of Hitchin's equations $(\delbar_E, t \varphi, h_t)$ \cite[Theorem~6.7]{MSWW14}. As shown in Fig.~\ref{fig:globalgluing}, the approximate metric $h_t^{\mathrm{app}}$ is constructed by desingularizing the singular metric $h_\infty$ by gluing in model solutions on disks around the zeros of~$q_2$. (These model solutions appear in Exercise~\ref{ex:model}.) Thus, they define the ``approximate Hitchin moduli space'' $\cM'_{\app}$ to be the moduli space of triples $\big[\big(\delbar_E, t \varphi, h_t^\app\big)\big]$. It too has a natural (non-hyperk\"ahler) $L^2$-metric~$g_{\app}$.

Mazzeo--Swoboda--Weiss--Witt decompose the difference $g_{\cM} - g_{\semif}$ into two pieces
 \begin{gather*}
 g_{\cM}-g_{\semif} = (g_{\cM} - g_{\app} ) + ( g_{\app} - g_{\semif}).
 \end{gather*}
They prove that $g_{\cM}-g_{\app}$ is exponentially decaying. All of their possible polynomial terms come from the second term, $g_{\app} - g_{\semif}$. Moreover, since $h_t^{\app}=h_\infty$ on the complement of the disks, the difference of the two metrics $g_{\app} - g_{\semif}$ reduces to an integral on disks around the ramification points.

\textbf{Idea 3.} Essentially\footnote{Dumas--Neitzke do not actually use the approximate solutions of \cite{MSWW14}.}, Dumas--Neitzke have a very clever way of dealing with the term \smash{$g_{\app}- g_{\semif}$} on the disks. The possible polynomial terms in Mazzeo--Swoboda--Weiss--Witt's asym\-p\-totic expansion are roughly from variations in which the zeros of $\det(\varphi + \eps \dot{\varphi})$ move. Roughly, Dumas--Neitzke use a local biholomorphic flow on the disks around each zero of $q_2$ that perfectly matches the changing location of the zero of $q_2+ \eps \dot{q}_2$.

The proof in \cite{Fredricksonasymptoticgeometry} can be seen as an extension of the method of Dumas--Neitzke~\cite{DumasNeitzke} to all of~$\cM'$ using the analysis and approximate solutions of Mazzeo--Swoboda--Weiss--Witt in \cite{MSWW14, MSWW17}.

\subsection*{Acknowledgements}
These notes are based on 3-hour mini-course aimed at early graduate students, given on November 11--12, 2017 at UIC. This course was part of the workshop ``Workshop on the geometry and physics of Higgs bundles'' and following conference ``Current Trends for Spectral Data III'' organized by Laura Schaposnik. (The notes have been updated to include a survey of results through October 2018.) My trip for the mini-course was funded by: the UIC NSF RTG grant DMS-1246844; L.P. Schaposnik's UIC Start up fund; and, NSF DMS 1107452, 1107263, 1107367 “RNMS: GEometric structures And Representation varieties“ (the GEAR Network).

I thank Laura Schaposnik for organizing the events and for her encouragement to contribute these notes. I thank Rafe Mazzeo for many discussions about the asymptotic geometry of the Hitchin moduli space, and Rafe Mazzeo, Steve Rayan, and the anonymous referees for their useful suggestions and comments.

\pdfbookmark[1]{References}{ref}
\LastPageEnding

\end{document}